\documentclass[a4paper,11pt]{article}
\usepackage[latin1]{inputenc}
\usepackage[english]{babel}
\usepackage{amsmath}
\usepackage{amsfonts}
\usepackage{amssymb}
\usepackage{epsfig}
\usepackage{amsopn}
\usepackage{amsthm}
\usepackage{color}
\usepackage{graphicx}
\newtheorem{theorem}{Theorem}
\newtheorem{lemma}{Lemma}
\newtheorem{remark}{Remark}

\newtheorem{proposition}{Proposition}
\newtheorem{corollary}{Corollary}
\newtheorem*{theorem*}{Theorem}
\newtheorem*{lemma*}{Lemma}
\newtheorem*{remark*}{Remark}
\newtheorem*{definition*}{Definition}
\newtheorem*{proposition*}{Proposition}
\newtheorem*{corollary*}{Corollary}
\numberwithin{equation}{section}
\numberwithin{theorem}{section}
\numberwithin{proposition}{section}
\numberwithin{lemma}{section}
\numberwithin{definition}{section}
\numberwithin{corollary}{section}
\numberwithin{remark}{section}

\newcommand{\real}{\mathbb{R}}

\def\a{\alpha}

\def\e{\varepsilon}        

\def\s{\sigma}          

\def\cp{{\cal P}}

\newcommand{\rd}{{\rm d}}

\newcommand{\RR}{\mathbb{R}}

\def\qed{\,\unskip\kern 6pt \penalty 500
\raise -2pt\hbox{\vrule \vbox to8pt{\hrule width 6pt
\vfill\hrule}\vrule}\par}
\definecolor{darkblue}{rgb}{0.05, .05, .65}
\definecolor{darkgreen}{rgb}{0.1, .65, .1}
\definecolor{darkred}{rgb}{0.8,0,0}
\newcommand{\beqn}{\begin{equation}}
\newcommand{\eeqn}{\end{equation}}
\newcommand{\bear}{\begin{eqnarray}}
\newcommand{\eear}{\end{eqnarray}}
\newcommand{\bean}{\begin{eqnarray*}}
\newcommand{\eean}{\end{eqnarray*}}
\begin{document}

\title{\huge \bf  Asymptotic behaviour of a nonlinear parabolic equation  with gradient absorption and critical exponent}

\author{
\Large Razvan Gabriel Iagar\,\footnote{Departamento de
Matem\'aticas, Universidad Aut\'onoma de Madrid, Campus de
Cantoblanco, E--28049 Madrid, Spain. \textit{e-mail:}
razvan.iagar@uam.es},\\[4pt] \Large Philippe Lauren\c cot\,\footnote{Institut de
Math\'ematiques de Toulouse, CNRS UMR~5219, Universit\'e de
Toulouse, F--31062 Toulouse Cedex 9, France. \textit{e-mail:}
Philippe.Laurencot@math.univ-toulouse.fr},\\ [4pt] \Large Juan Luis
V{\'a}zquez\,\footnote{Departamento de Matem\'aticas, Universidad
Aut\'onoma de Madrid, Campus de Cantoblanco, E--28049 Madrid, Spain.
Also affiliated with ICMAT, Madrid. \textit{e-mail:}
juanluis.vazquez@uam.es} }
\date{}
\maketitle

\begin{abstract}
We study the large-time behaviour of the solutions of the evolution
equation involving nonlinear diffusion and gradient absorption,
$$
\partial_t u - \Delta_p u + |\nabla u|^q=0\,.
$$
We consider the problem posed for $x\in \real^N $ and $t>0$ with
nonnegative and compactly supported initial data. We take the
exponent $p>2$ which corresponds to slow $p$-Laplacian diffusion.
The main feature of the paper is that the exponent $q$ takes the
critical value $q=p-1$ which leads to interesting asymptotics. This
is due to the fact that in this case both the Hamilton-Jacobi term $
|\nabla u|^q$ and the diffusive term $\Delta_p u$ have a similar
size for large times. The study performed in this paper shows that a
delicate asymptotic equilibrium happens, so that the large-time
behaviour of the solutions is described by a rescaled version of a
suitable self-similar solution of the Hamilton-Jacobi equation
$|\nabla W|^{p-1}=W$, with  logarithmic time corrections. The
asymptotic rescaled profile is a kind of sandpile with a cusp on
top, and it is independent of the space dimension.
\end{abstract}

\vspace{2.0 cm}


\medskip


\newpage

\section{Introduction and main results}\label{sect.intro}

In this paper we deal with the Cauchy problem associated to the
diffusion-absorption equation:
\begin{equation}  \label{a1}
\partial_t u - \Delta_p u + |\nabla u|^{q} = 0\,, \quad (t,x)\in Q\,,
\end{equation}
posed in  $Q:= (0,\infty)\times\RR^N$ with initial data
\begin{equation}\label{a2}
 u(0,x) = u_0(x)\ge 0\,,  \quad x\in\RR^N\,,
\end{equation}
where the $p$-Laplacian operator is
defined as usual by \ $
\Delta_p u := \mbox{ div }\left( |\nabla u|^{p-2}\ \nabla u \right).
$
To be specific we take $p>2$, which implies finite speed of
propagation, and we consider nonnegative weak solutions $u\ge 0$
with compactly supported initial data  $u_0$ such that
\begin{equation} \label{a3} u_0\in W^{1,\infty}(\RR^N)\,, \;\;
u_0\ge 0\,, \;\; \mbox{ supp }(u_0) \subset B(0,R_0)\,, \;\;
u_0\not\equiv 0\,,
\end{equation} for some $R_0>0$. Known properties of the equation ensure
that the corresponding solutions will be compactly supported with
respect to the space variable for every time $t>0$. The goal of the
paper is to describe in detail the asymptotic behaviour of the
solutions as $t\to\infty$.

The equation \eqref{a1} has been studied by various authors for different
values of the parameters $p\ge 2$ and $q>1$ as a model of linear
or nonlinear diffusion with gradient-dependent absorption, see
\cite{BKL04,BGK04,GL07,Gi05,La08} for the semilinear case $p=2$, and
\cite{ATU04,BtL08,LV07,Sh04a} for the quasilinear case $p>2$.  It has
been shown that the large-time behaviour of this initial-value
problem depends on the relative influence of the diffusion and
absorption terms and leads to a classification into the following
ranges of $q$:

\noindent (i) when $q>q_2:=p-N/(N+1)$ the large time behaviour is
purely diffusive and the first-order absorption term disappears in
the limit $t\to\infty$; this is a case of asymptotic simplification
in the sense of \cite{V91}.

\noindent (ii) For $q_1:=p-1<q<q_2$ there is a behaviour given by a
certain balance of diffusion and absorption in the form of a
self-similar solution, its existence being established in
\cite{Sh04a}; there is no asymptotic simplification;

\noindent (iii) for $1<q<p-1$ the two last authors have recently
shown in \cite{LV07} that the main term is the absorption term,
leading to a separate-variables asymptotic behaviour, with diffusion
playing a secondary role. We thus have asymptotic simplification, now
with absorption as the dominating effect.

The two critical cases $q=q_2$ and $q=q_1$ represent limit
behaviours, and as is often the case in such situations, they give
rise to interesting dynamics due to the curious interaction of two
effects of similar power. Such situations  usually lead to phenomena
called \textsl{resonances} in mechanics, with interesting
non-trivial mathematical analysis. Such interesting behaviour has
been shown in particular in \cite{GL07} for $q=q_2$, in the linear
case $p=2$, with the result that logarithmic factors modify the
purely diffusive behaviour found for $q>q_2$. A similar situation is
expected to be met when $p>2$ and $q=q_2$.

We devote this paper to describe the other limit case $q=q_1=p-1$
when $p>2$, the latter condition guaranteeing that $q>1$. In that
case the diffusion and the first order term have similar asymptotic
size and logarithmic corrections appear in the asymptotic rates. The
mathematical analysis that we  perform below is strongly tied to a
good knowledge of the expansion of the support of the solution, or
in other words, the location of the free boundary, which happens to
be approximately a sphere of radius $|x|\sim C\log t$ for large
times $t$. From now on, we assume that
$$
q =q_1=p-1\,.
$$

\subsection{Bounds in suitable norms}

Studying the large time behaviour of solutions and interfaces of our
problem relies on suitable and very precise estimates. The time
expansion of the support and the time decay of solutions to the
Cauchy problem \eqref{a1}-\eqref{a2}, with non-negative and
compactly supported initial data have been recently investigated in
\cite{BtL08}. The following results are proved:

\begin{proposition}\label{pra1}
Under the above assumptions on the equation and data, the Cauchy
problem \eqref{a1} has a unique non-negative viscosity solution
$$
u\in\mathcal{BC}([0,\infty)\times\RR^N)\cap
L^\infty(0,\infty;W^{1,\infty}(\RR^N))
$$
which satisfies:
\begin{equation} \label{a4} 0 \le u(t,x) \le \|u_0\|_\infty\,, \quad
(t,x)\in Q\,,
\end{equation}
\begin{equation}\label{a4b}
\|\nabla u(t)\|_\infty \le \|\nabla u_0\|_\infty\,, \quad t\ge 0\,,
\end{equation}
\begin{equation} \label{a5} \mbox{ supp
}(u(t)) \subset B(0,C_1\ \log{t}) \;\;\mbox{ for all }\;\; t\ge 2\,,
\end{equation}
together with the following norm estimates
\bear \label{a6}
 \|u(t)\|_1 & \le & C_2\ t^{-1/(p-2)}\ (\log{t})^{(p(N+1)-2N-1)/(p-2)} \;\;\mbox{ for all }\;\; t\ge 2\,, \\
 \label{a7}
 \|u(t)\|_\infty & \le & C_2\ t^{-1/(p-2)}\ (\log{t})^{(p-1)/(p-2)} \;\;\mbox{ for all }\;\; t\ge 2\,, \\
 \label{a8}
 \|\nabla u(t)\|_\infty & \le & C_2\ t^{-1/(p-2)}\ (\log{t})^{1/(p-2)} \;\;\mbox{ for all }\;\; t\ge 2\,,
\eear for some positive constants $C_1$ and $C_2$ depending only on
$p$, $N$, and $u_0$.
\end{proposition}
Here and below, $\mathcal{BC}([0,\infty)\times\RR^N)$ denotes the
space of bounded and continuous functions on $[0,\infty)\times\RR^N$
and $\Vert \cdot\Vert_r$ denotes the $L^r(\RR^N)$-norm for $r\in
[1,\infty]$. As we shall see, these bounds will be very useful in the
sequel. The well-posedness of \eqref{a1}-\eqref{a2} and the
properties \eqref{a4}, \eqref{a5}, and \eqref{a6} are established in
\cite[Theorems~1.1 \&~1.6, Corollary~1.7]{BtL08}, while \eqref{a7}
and \eqref{a8} follow from \eqref{a6} and
\cite[Proposition~1.4]{BtL08}. We will also use the notation
$r_{+}=\max\{r,0\}$ for the positive part of the real number $r$.

\subsection{Main results}\label{Sect.main}

We describe next the main contribution of this paper. As already
mentioned, our goal is to study the asymptotic behaviour of the
solution $u$ of the resonant problem \eqref{a1} with $p>2$ and
$q=p-1$, and with compactly supported and nonnegative initial data.
Moreover, since the equation has the property of finite
speed of propagation, it is natural to raise the question about how
the interface and the support of the solution expand in time. We
also answer this question in the present paper.

\medskip

\noindent \textbf{Asymptotic behaviour.} The main result is the
following:
\begin{theorem}\label{asympt.main}
Let $u$ be the solution of the Cauchy problem \eqref{a1}-\eqref{a2},
with $u_0$ as in \eqref{a3}. Then, $u$ decays in time like \
$O(t^{-1/(p-2)}(\log t)^{(p-1)/(p-2)})$  \ and the support spreads in space like
\ $O(\log t)$ as $t\to\infty$. More precisely, we have the limit:
\begin{equation}\label{main.asympt}
\lim_{t\to\infty}\sup_{x\in\real^N}\left|
\frac{c_p\, t^{1/(p-2)}}{(\log{t})^{(p-1)/(p-2)}}\ u\left(t,x\right) -
\ \left(1-\frac{(p-2)|x|}{\log\,t}\right)_{+}^{(p-1)/(p-2)} \right| = 0\,,
\end{equation}
with precise constant
\begin{equation*}
c_p=(p-2)^{1/(p-2)}(p-1)^{(p-1)/(p-2)}.
\end{equation*}
\end{theorem}

In the proof, the expression of the asymptotic profile is obtained after
a complicated time scaling of $u$ and $x$ in the form of uniform limit
\begin{equation}\label{main.asympt2}
\frac{t^{1/(p-2)}}{(\log{t})^{(p-1)/(p-2)}}\ u\left(t,x \right) \to
(p-2)^{-p/(p-2)}\ W((p-2)x/\log t),
\end{equation}
where the function
\begin{equation}
W(x):=\left(\frac{p-2}{p-1}\ (1-|x|)_{+}\right)^{(p-1)/(p-2)}
\end{equation}
is the unique viscosity solution to the stationary form of the rescaled problem, which is:
\begin{equation}\label{eqlim}
|\nabla W|^{p-1}-W=0 \ \mbox{ in } \ B(0,1)\,, \quad W=0 \ \mbox{ on
} \ \partial B(0,1), \ \quad W>0 \ \mbox{ in } \ B(0,1).
\end{equation}

Let us notice that, as usual in resonance cases, the limit profile
is not a self-similar solution, but it introduces logarithmic
corrections to a self-similar, separate-variables profile (which in
our case is $t^{-1/(p-2)}(p-2)^{-p/(p-2)}W((p-2)x)$). The uniqueness
of $W$ as viscosity solution of \eqref{eqlim} is very important in
the proof and follows from \cite{I87}.

In consonance with \eqref{main.asympt}, we show that the shape
of the support of $u(t)$ gets closer to a ball while expanding as
time goes by. This is in sharp contrast with the situation described
in \cite{LV07} for \eqref{a1} in the intermediate range $q\in
(1,p-1)$, $p>2$ where the positivity set stays bounded and can have
a very general shape. When $q=p-1$, the diffusion thus acts in three
directions: the scaling is different, the support grows unboundedly
with time, and the geometry of the positivity set simplifies. Another remarkable
consequence of the interplay diffusion-absorption is the fact that the asymptotic
profile is radially symmetric and does not depend on the space dimension.

We devote Section~\ref{Sect.scaling2} to the proof of
Theorem~\ref{asympt.main}. For the proof, we use a precise estimate
for the propagation of the positivity set, that is described below.
Another tool is the existence of a large family of subsolutions
having a special, explicit form and allowing for a theoretical argument with viscosity
solutions to finish the proof.

\medskip

\noindent \textbf{Propagation of the positivity set.} We denote the
positivity set and its maximal expansion radius by
\begin{equation}
\cp_{u}(t):=\{x\in\real^N: \ u(t,x)>0\}, \qquad \gamma(t)=\sup\{|x|:
\ x\in P_{u}(t)\}
\end{equation}
respectively. Then:
\begin{theorem}\label{main.posit}
Under the running notations and assumptions, we have:
\begin{equation}
\lim\limits_{t\to\infty}\frac{\gamma(t)}{\log\,t}=\frac{1}{p-2}.
\end{equation}
Moreover, the free boundary of $u$ has the same speed of expansion
in any given direction $\omega\in\real^N$ with $|\omega|=1$.
\end{theorem}
In fact, we give more precise estimates for the expansion of the
positivity region, obtained via comparison with some well-chosen
traveling waves. The proof of Theorem~\ref{main.posit} is performed
in Section~\ref{sect.posit}.

\medskip

\noindent \textbf{Two scalings.} In order to prove the two theorems, we have to perform two different
scaling steps. The first scaling, described in formula (\ref{a9}) below, is
the natural one corresponding to standard scaling invariance; such a scaling
 has also been used in \cite{LV07} in the case $q\in(1,p-1)$ to obtain the
correct scale of the solutions. But for $q=p-1$, we observe that a
phenomenon of grow-up appears, which is typical for resonance cases:
the effect of the resonance implies that the rescaled solution does
not stabilize in time; on the contrary, it grows and becomes
unbounded in infinite time. That is why we need a second scaling,
given by the new functions $w$ and $y$ defined in (\ref{def.y}) and
(\ref{def.w}), which is less natural but turns out to be adapted
to our problem: it takes into account the logarithmic corrections
(suggested by the a priori estimates of Proposition~\ref{pra1},
which turn out to be sharp), and it is adapted to the size of the
grow-up phenomenon; thus, in the rescaled variables we can describe
the real form and behaviour of the solution.


\section{Scaling variables I}\label{Sect.Scaling1}

We recall that $p>2$ and $q=p-1$. We introduce a first set of
\textit{self-similar} variables; we keep the space variable $x$
and introduce logarithmic time
\begin{equation}
\tau:=\frac{1}{p-2}\,\ln{(1+(p-2)t)},
\end{equation}
as well as the new unknown function $v=v(\tau,x)$ defined by
\begin{equation} \label{a9}
u(t,x) = (1+(p-2)t)^{-1/(p-2)}\ v\left( \tau, x \right)\,, \quad
(t,x)\in [0,\infty)\times\RR^N\,.
\end{equation}
Clearly, $v$ solves the rescaled equation
\begin{equation}
\partial_\tau v - \Delta_p v + |\nabla v|^q -v  =  0\ \,,
\quad (\tau,x)\in Q\,, \label{a10}
\end{equation}
with the same initial condition
\begin{equation}
\label{a11} v(0) =  u_0\,, \quad x\in\RR^N\,.
\end{equation}
We next translate the a priori bounds  (\ref{a6}), (\ref{a7}), and
(\ref{a8}) in terms of the rescaled function $v$: there is $C_3>0$
depending only on $p$, $N$, and $u_0$ such that
\begin{equation}
\label{a12} \frac{\Vert v(\tau)\Vert_1}{\tau^{(p(N+1)-2N-1)/(p-2)}}
+ \frac{\Vert v(\tau)\Vert_\infty}{\tau^{(p-1)/(p-2)}} + \frac{\Vert
\nabla v(\tau)\Vert_\infty}{\tau^{1/(p-2)}} \le C_3 \;\;\mbox{ for
}\;\; \tau\ge 1\,.
\end{equation}

\subsection{The positivity set: time monotonicity}

We define the positivity set $\mathcal{P}_v(\tau)$ of the function
$v$ at time $\tau\ge 0$ by
\begin{equation}\label{b1}
\mathcal{P}_v(\tau):= \left\{ x\in\RR^N : \;\;
v(\tau,x)>0  \right\}\,.
\end{equation}
\begin{proposition}\label{prb1}
For $\tau_1\in [0,\infty)$ and $\tau_2\in (\tau_1,\infty)$ we have
\begin{equation}\label{b2}
\mathcal{P}_v(\tau_1) \subseteq \mathcal{P}_v(\tau_2)
\quad\mbox{ and }\quad \bigcup_{\tau\ge 0}
\mathcal{P}_v(\tau)=\RR^N\,.
\end{equation}
In addition, for each $x\in\RR^N$ there are $T_x\ge 0$ and
$\varepsilon_x>0$ such that
\begin{equation}\label{b3}
v(\tau,x) \ge \varepsilon_x\ \tau^{(p-1)/(p-2)} \;\;\mbox{ for }\;\;
\tau\ge T_x\,.
\end{equation}
\end{proposition}
The proof relies on the availability of suitable subsolutions which
we describe next.

\begin{lemma}\label{leb2}
Define two positive real numbers $R_p$ and $T_p$ by
$$
R_p := \frac{p-2}{2^p (p-1)} \quad\mbox{ and }\quad T_p :=
\frac{2(p-1)}{p-2}\ \left( 2 + 2^{p-1} (N+p-2) \right)\,.
$$
If $R\in (0,R_p]$ and $T\ge T_p$, the function $s_{R,T}$ given by
$$
s_{R,T}(\tau,x) := \frac{p-2}{R (p-1)}\ (T+\tau)^{(p-1)/(p-2)}\
\left( R^2 - \frac{|x|^2}{(T+\tau)^2} \right)_+^{(p-1)/(p-2)}, \quad
(\tau,x)\in Q\,,
$$
is a (viscosity) subsolution to \eqref{a10}.
\end{lemma}

\noindent\textbf{Proof.} We have $s_{R,T}(\tau,x)=
(T+\tau)^{(p-1)/(p-2)}\ \sigma(\xi)$ with $\xi:=x/(T+\tau)$ and
$\sigma(\xi):= (p-2)\ \left( R^2 - |\xi|^2 \right)_+^{(p-1)/(p-2)}
/(R(p-1))$. Since $p-1>p-2>0$, we observe that $\sigma$ and
$|\nabla\sigma|^{p-2} \nabla\sigma$ both belong to
$\mathcal{C}^1(\RR^N)$. Therefore,
$$
L(\tau,x) := R\ (T+\tau)^{-(p-1)/(p-2)}\ \left\{ \partial_\tau
s_{R,T} - \Delta_p s_{R,T} + \left| \nabla s_{R,T}\right|^{p-1} -
s_{R,T} \right\}
$$
is well-defined for $(\tau,x)\in [0,\infty)\times\RR^N$ and \bean
L(\tau,x) & = & \frac{R}{T+\tau}\ \left\{ \frac{p-1}{p-2}\ \sigma(\xi) - \xi\cdot\nabla\sigma(\xi) - \Delta_p \sigma(\xi) \right\} + R\ \left| \nabla \sigma(\xi)\right|^{p-1} - R\ \sigma(\xi) \\
& = & \left( R^2 - |\xi|^2 \right)_+^{(p-1)/(p-2)}\ \left\{ \frac{1}{T+\tau}\ \left( 1 + 2^{p-1} (N+p-2)\ \frac{|\xi|^{p-2}}{R^{p-2}} \right) \right\} \\
& + & \left( R^2 - |\xi|^2 \right)_+^{(p-1)/(p-2)}\ \left\{ \frac{2}{T+\tau}\ \frac{|\xi|^2}{R^2-|\xi|^2}\ \left( 1 - \frac{2^{p-1} (p-1)}{p-2}\ \frac{|\xi|^{p-2}}{R^{p-2}} \right) \right\} \\
& + & \left( R^2 - |\xi|^2 \right)_+^{(p-1)/(p-2)}\ \left\{ 2^{p-1} \frac{|\xi|^{p-1}}{R^{p-2}} - \frac{p-2}{p-1} \right\} \\
& \le & \left( R^2 - |\xi|^2 \right)_+^{(p-1)/(p-2)}\ \left\{ \frac{1+2^{p-1} (N+p-2)}{T} + 2^{p-1} R - \frac{p-2}{p-1} \right\} \\
& + & \left( R^2 - |\xi|^2 \right)_+^{(p-1)/(p-2)}\ \left\{
\frac{2}{T+\tau}\ \frac{|\xi|^2}{R^2-|\xi|^2}\ \left( 1 -
\frac{2^{p-1} (p-1)}{p-2}\ \frac{|\xi|^{p-2}}{R^{p-2}} \right)_+
\right\} \,. \eean
We next note that
$$
1 - \frac{2^{p-1} (p-1)}{p-2}\ \frac{|\xi|^{p-2}}{R^{p-2}}\le 0
\quad\mbox{ if }\quad |\xi|\ge \frac{R}{2}\,,
$$
so that the last term of the right-hand side of the previous
inequality is bounded from above by $2\left( R^2 - |\xi|^2
\right)_+^{(p-1)/(p-2)}/(3T)$. Consequently, owing to the choice of
$R$ and $T$, \bean
L(\tau,x) & \le & \left( R^2 - |\xi|^2 \right)_+^{(p-1)/(p-2)}\ \left\{ \frac{1+2^{p-1} (N+p-2)}{T_p} + 2^{p-1} R_p - \frac{p-2}{p-1} + \frac{2}{3T_p} \right\} \\
& \le & 0\,, \eean whence the claim. \qed

\medskip

\noindent\textbf{Proof of Proposition~\ref{prb1}.} (i) Fix $\tau_1\ge
0$ and $x_1\in\mathcal{P}_v(\tau_1)$. Owing to the continuity of
$x\longmapsto v(\tau_1,x)$ there are $\delta>0$ and $r_1>0$ such
that $v(\tau_1,x)\ge\delta$ for $x\in B(x_1,r_1)$. Take now $R>0$
small enough such that $R<\min{\{r_1,R_p\}}$ and satisfying
$$
R< \frac{r_1}{T_p+\tau_1} \quad\mbox{ and }\quad \frac{p-2}{p-1}\
(T_p+\tau_1)^{(p-1)/(p-2)}\ R^{p/(p-2)} \le \delta\,,
$$
the parameters $R_p$ and $T_p$ being defined in Lemma~\ref{leb2}.
Then we have $s_{R,T_p}(\tau_1,x-x_1) = 0 \le v(\tau_1,x)$ \ if \ $|x-x_1|\ge R\ (T_p+\tau_1)$, while
$$
s_{R,T_p}(\tau_1,x-x_1) \le \frac{p-2}{R(p-1)}\
(T_p+\tau_1)^{(p-1)/(p-2)}\ R^{(2p-2)/(p-2)} \le \delta \le
v(\tau_1,x)
$$
if $|x-x_1|\le R\ (T_p+\tau_1)$ as $R(T_p+\tau_1)\le r_1$. Moreover,
if $\tau_2>\tau_1$,  $\tau\in [\tau_1,\tau_2]$ and $x\in\partial
B(x_1,R (T_p+\tau_2))$, then $s_{R,T_p}(\tau,x-x_1)=0 \le
v(\tau,x)$. Recalling that $s_{R,T_p}$ is a subsolution to
\eqref{a10} by Lemma~\ref{leb2}, we infer from the comparison
principle that $s_{R,T_p}(\tau,x-x_1) \le v(\tau,x)$ for
$(\tau,x)\in [\tau_1,\tau_2]\times B(x_1,R (T_p+\tau_2))$. As
$s_{R,T_p}(\tau,x-x_1)=0\le v(\tau,x)$ for $\tau\in [\tau_1,\tau_2]$
and $x\not\in B(x_1,R (T_p+\tau_2))$ we actually have
$s_{R,T_p}(\tau,x-x_1) \le v(\tau,x)$ for $(\tau,x)\in
[\tau_1,\tau_2]\times\RR^N$. Since $\tau_2>\tau_1$ is arbitrary and
neither $R$ nor $T_p$ depend on $\tau_2$, we end up with
\begin{equation}\label{bs31}
s_{R,T_p}(\tau,x-x_1) \le v(\tau,x)\,, \qquad (\tau,x)\in
[\tau_1,\infty)\times\RR^N\,.
\end{equation}
A first consequence of \eqref{bs31} is that, if $\tau_2>\tau_1$,
then $v(\tau_2,x_1)\ge s_{R,T_p}(\tau_2,0)>0$ so that $x_1$ also
belongs to  $\mathcal{P}_v(\tau_2)$.

Next, given $x\in\RR^N$, we have $x\in B(x_1,R(T_p+\tau))$ for
$\tau$ large enough  and it follows from \eqref{bs31} that
$v(\tau,x)\ge s_{R,T_p}(\tau,x-x_1) > 0$ for $\tau$
large enough. Consequently, $x$ belongs to $\mathcal{P}_v(\tau)$ for
$\tau$ large enough which proves the second assertion of \eqref{b2}.

\medskip

\noindent (ii)  Consider $x_0\in\RR^N$. According to \eqref{b2} there is
$\tau_0$ large enough such that $x_0\in\mathcal{P}_v(\tau_0)$.
Arguing as in the proof of \eqref{b2}, we may find $r_0$ small
enough (depending on $x_0$) such that $s_{r_0,T_p}(\tau,x-x_0) \le
v(\tau,x)$ for $(\tau,x)\in [\tau_0,\infty)\times\RR^N$.
Consequently,
$$
v(\tau,x_0) \ge \frac{p-2}{r_0 (p-1)}\ (T_p+\tau)^{(p-1)/(p-2)}\
r_0^{(2p-2)/(p-2)} \ge \frac{p-2}{p-1}\ r_0^{p/(p-2)}\
\tau^{(p-1)/(p-2)}\,,
$$
which gives the lower bound \eqref{b3}.  \qed

\medskip

\begin{corollary}\label{cor:bfbsb}
Assume that $u_0(0)>0$. Then there is $r_*>0$ such that
\begin{equation}\label{spirou}
v(\tau,x) \ge \frac{(p-2)}{r_* (p-1)} (1+\tau)^{(p-1)/(p-2)} \left( r_*^2 - \frac{|x|^2}{(1+\tau)^2} \right)_+^{(p-1)/(p-2)}\,, \quad (\tau,x)\in Q\,.
\end{equation}
\end{corollary}

\noindent\textbf{Proof.} Arguing as in the proof of \eqref{b2} and
using the positivity of $u_0(0)$, we may find $r_*>0$ small enough
such that $s_{r_*,T_p}(\tau,x)\le v(\tau,x)$ for $(\tau,x)\in Q$.
Since $T_p>1$, the previous inequality implies \eqref{spirou}. \qed

\subsection{Eventual radial symmetry}

We prove the following classical monotonicity lemma, see
\cite[Proposition~2.1]{AC83} for instance.
\begin{lemma}\label{leb4}
If $x\in\RR^N$ and $r>0$ satisfy \ $ |x|>2R_0$ and $ r< |x| - 2 R_0$. Then,
\begin{equation}
v(\tau,x) \le \inf_{|y|=r} v(\tau,y)  \quad\mbox{ for }\quad \tau\ge
0\,.
\end{equation}
Here, $R_0$ is radius of the initial ball defined in \eqref{a3}.
\end{lemma}

\noindent\textbf{Proof.} The proof relies on Alexandrov's
reflection principle. Let $(x,r)\in\RR^N\times (0,\infty)$ fulfil
the assumptions of Lemma~\ref{leb4} and consider $y\in\RR^N$ such
that $|y|=r$. Let $H$ be the hyperplane of points of $\RR^N$ which
are equidistant from $x$ and $y$, namely
$$
H := \left\{ z\in\RR^N\ : \ \left\langle z - \frac{x+y}{2} , x-y
\right\rangle = 0 \right\}\,.
$$
Introducing
$$
H_- := \left\{ z\in\RR^N\ : \ \left\langle z - \frac{x+y}{2} , x-y
\right\rangle \le 0 \right\}
$$
and
$$
\tilde{v}(\tau,z) := v\left( \tau , z - 2\ \left\langle z -
\frac{x+y}{2} , x-y \right\rangle\ \frac{x-y}{|x-y|^2} \right)\,,
\qquad (\tau,z)\in Q\,,
$$
it readily follows from the rotational and translational invariance
of \eqref{a10} that $\tilde{v}$ also solves \eqref{a10}. In
addition, $y\in H_-$ and $\mathcal{P}_v(0) \subseteq B(0,R_0)
\subseteq H_-$ by \eqref{a3}. Now, on the one hand, if $z\in H_-$,
then
$$
z - 2\ \left\langle z - \frac{x+y}{2} , x-y \right\rangle\
\frac{x-y}{|x-y|^2} \not\in H_-
$$
and $\tilde{v}(0,z)=0\le v(0,z)$. On the other hand, if $z\in H =
\partial H_-$ and $\tau\ge 0$, we clearly have
$\tilde{v}(\tau,z)=v(\tau,z)$. We are then in a position to apply
the comparison principle to \eqref{a10} on $(0,\infty)\times H_-$
and conclude that
\begin{equation}\label{b5}
\tilde{v}(\tau,z) \le v(\tau,z)\,, \qquad (\tau,z)\in
[0,\infty)\times H_-\,.
\end{equation}
Recalling that $y\in H_-$, we infer from \eqref{b5} that
$v(\tau,y)\ge\tilde{v}(\tau,y)=v(\tau,x)$ for $\tau\ge 0$ which is
the expected result. \qed

\medskip

\begin{remark}
Although Lemma~\ref{leb4} will not be used in the main proofs, this
is an interesting result for the qualitative theory, since it shows
that the dynamics symmetrizes the solution.
\end{remark}

\section{Propagation of the positivity set}\label{sect.posit}

We next turn to the speed of expansion of the positivity set
$\mathcal{P}_v$ of $v$ and put
\begin{equation}\label{bs32}
\varrho(\tau) := \sup{\left\{ |x|\ : \ x \in \mathcal{P}_v(\tau)
\right\}}\,,
\end{equation}
so that $\mathcal{P}_v(\tau)\subseteq B(0,\varrho(\tau))$ for
$\tau\ge 0$. The purpose of this section is to prove that the
expansion speed $\varrho(\tau)$ of $\mathcal{P}_v(\tau)$ is
asymptotically equal to $\tau$, in other words,
\begin{equation*}
\lim\limits_{\tau\to\infty}\frac{\varrho(\tau)}{\tau}=1,
\end{equation*}
and, more precisely, to prove Theorem \ref{main.posit}.

The proof relies on the existence of ``nice'' traveling wave
solutions of \eqref{a10}, which may be used as subsolutions and
supersolutions for the Cauchy problem \eqref{a10}-\eqref{a11}. The
construction of such traveling waves is inspired on the technique
used in the so-called KPP problems, \cite{KPP}, which has developed
a wide literature; see e.\,g., \cite{Ar80}, \cite{VaTube} for
applications to porous media, and \cite{QRV} for blow-up problems.
We thus begin with a phase-plane analysis, proving the existence of
the desired traveling waves.

\subsection{Traveling wave analysis for $N=1$}\label{sect.TW}

We look for traveling waves of the form
\begin{equation*}
v(\tau,x)=f(z), \quad z=x-c\tau, \ c>0,
\end{equation*}
solving \eqref{a10} in dimension $N=1$. Then, the profile $f$ solves
the ordinary differential equation:
\begin{equation}\label{OdeTW}
-cf'-\left(|f'|^{p-2}f'\right)'+|f'|^{p-1}-f=0.
\end{equation}
We are actually only interested  in traveling waves which present an
interface, that is, $f$ vanishes for $z$ sufficiently large. As we
shall see below, the profile $f$ is non-monotone in general, but is
nonnegative and decreasing near the interface. We transform
\eqref{OdeTW} into a first order system, by introducing the notation
$U=f$ and $V=-f'$. We arrive at the following system
\begin{equation}\label{syst1}
\left\{\begin{array}{ll}(p-1)|V|^{p-2}U'=-(p-1)|V|^{p-2}V, \\
(p-1)|V|^{p-2}V'=-cV-|V|^{p-1}+U,
\end{array}\right.
\end{equation}
where, for the orbits, the term $(p-1)|V|^{p-2}$ in the right-hand
side has no influence (since we work with $dV/dU$) and can be
ignored after a change of the time variable. We perform next the
phase-plane analysis of the system \eqref{syst1}.

\medskip

\noindent \textbf{Local analysis in the plane}. The system
\eqref{syst1} has a unique critical point, $P=(0,0)$, and the
Jacobian matrix $J(0,0)$ at this point is given by
\begin{equation*}
J(0,0)=\left(
          \begin{array}{cc}
            0 & 0 \\
            1 & -c \\
          \end{array}
        \right)
\end{equation*}
with eigenvalues $\lambda_1=0$ and $\lambda_2=-c$, and corresponding
eigenvectors are $e_1=(c,1)$ and $e_2=(0,1)$.  By a careful
analysis, we notice that the center manifold in $P$ is tangent to
$e_1$, and is asymptotically stable. It follows that $P$ is a stable
node for every $c>0$. There is a unique orbit entering $P$ and
tangent to $e_2$, forming the stable manifold; its local behaviour
is $U(z)\sim C(-z)^{(p-1)/(p-2)}$ as $z\to 0$, hence this orbit
contains all the traveling waves with velocity $c$ and having an
interface. By standard theory (see, e.g., \cite{Pe}), all the other
orbits approach the center manifold, tangent to $e_1$, and present
an exponential decay, but no interface: $U(z)\sim e^{-cz}$ as
$z\to\infty$.

\medskip

\noindent \textbf{Local analysis at infinity}. We investigate the
behaviour of the system when $U$ is very large. For monotone
traveling waves, we make the following inversion of the plane:
\begin{equation*}
Z=\frac{1}{U}, \quad W=\frac{|V|^{p-2}V}{U},
\end{equation*}
and we are interested in the local behaviour near $Z=0$. After straightforward calculations, \eqref{syst1} becomes the new system:
\begin{equation}\label{inftysyst}
\left\{\begin{array}{ll} Z'=Z^{(2p-3)/(p-1)}W|W|^{-(p-2)/(p-1)},\\
W'=Z^{(p-2)/(p-1)}|W|^{p/(p-1)}-cZ^{(p-2)/(p-1)}W|W|^{-(p-2)/(p-1)}+1-|W|.\end{array}\right.
\end{equation}
We find two critical points with $Z=0$, namely $Q_1=(0,1)$ and
$Q_2=(0,-1)$. We will analyze only $Q_1$, i.e. the decreasing
traveling waves. Let us also remark that, in the second equation of
\eqref{inftysyst}, the terms with $Z$ are dominated by $1-|W|$ near
$Q_1$ and $Q_2$, hence we can study the local behaviour by using the
approximate equation only with $1-|W|$ in the right-hand side. The
linearization near $Q_1$ has eigenvalues $\lambda_1=0$ and
$\lambda_2=-1$, and the center manifold, which is tangent to the
line $W=1$, is unstable. Hence, the point $Q_1$ behaves like a
saddle, and the orbits which are interesting for our study are the
orbits going out of $Q_1$. These orbits are tangent to $W=1$, and in
the original system they satisfy $U\sim V^{p-1}$, hence, by
integration,
\begin{equation*}
U(z)\sim |z|^{(p-1)/(p-2)}, \quad \mbox{as} \ z\to-\infty,
\end{equation*}
and are decreasing. The local analysis around $Q_2$ is similar, but not interesting for our goals.

Let us notice that not all solutions passing through a point in the
first quadrant come from $Q_1$. Indeed, the orbits touching the
curve $U=cV+V^{p-1}$ change monotonicity as functions $V=V(U)$,
hence they have previously reached the axis $V=0$, meaning a change
of monotonicity as $f=f(z)$, and they enter through this change in
the first quadrant. Analyzing the curve $U=cV+V^{p-1}$, we observe
that it connects in the phase-plane the points $P=(0,0)$ and $Q_1$,
being tangent in $Q_1$ to the axis $W=1$. In particular, there exist
non-monotone solutions, and this is the object we are interested in.

\medskip

\noindent \textbf{Global behaviour}. This is now not difficult to
establish, by merging the previous local analysis with the following
important remarks:

\noindent (a) The evolution of the system \eqref{syst1} with respect
to the parameter $c$ is monotone. Indeed, we calculate:
\begin{equation*}
\frac{\rd}{\rd c}\left(\frac{\rd V}{\rd U}\right)=\frac{1}{(p-1)|V|^{p-2}}>0.
\end{equation*}

\noindent (b) There exists an explicit family of traveling wave solutions with speed $c=1$:
\begin{equation}\label{expl.TW}
f_{1,K}(z)=\left(\frac{p-2}{p-1}\right)^{(p-1)/(p-2)}(K-z)_{+}^{(p-1)/(p-2)}\,.
\quad K\ge 0,
\end{equation}
This function is obviously decreasing and presents an interface at
$z=K$. It is immediate to check that this orbit satisfies
$U=V^{p-1}$, hence it comes from the point $Q_1$ along the center
manifold of it, and it enters $P$, being the unique orbit entering
$P$ and tangent to the eigenvector $e_2=(0,1)$ (unique for $c=1$),
as discussed above.

\noindent (c) Moreover, the vectors of the direction field of
\eqref{syst1} over the curve $U=V^{p-1}$ (which gives the explicit
orbit \eqref{expl.TW}) have the same direction. Indeed, the normal
vector to this curve is $(1,-(p-1)V^{p-2})$ and we calculate:
\begin{equation*}
(1,-(p-1)V^{p-2})\cdot(-(p-1)V^{p-1},-cV-V^{p-1}+U)=(p-1)(c-1)V^{p-1}.
\end{equation*}
For $c=1$ we obtain the explicit trajectory, and for $c<1$, the
above scalar product is negative, hence all these vectors have the
same direction, contrary to $(1,-(p-1)V^{p-2})$. For $c>1$, all
these vectors have the same direction as $V$.

Since we are interested only in traveling waves with an interface,
we analyze only the unique (for $c$ fixed) orbit entering $P=(0,0)$
tangent to $e_2=(0,1)$. For $c=1$, it is explicit and connects $P$
and $Q_1$ in the first quadrant. We draw the phase-plane for $c=1$
in Figure~\ref{figure1} below; it is clear that the explicit
connection will not change sign and monotonicity.

\medskip

\begin{figure}[ht!]
   \begin{center}
   \includegraphics[width=8cm]{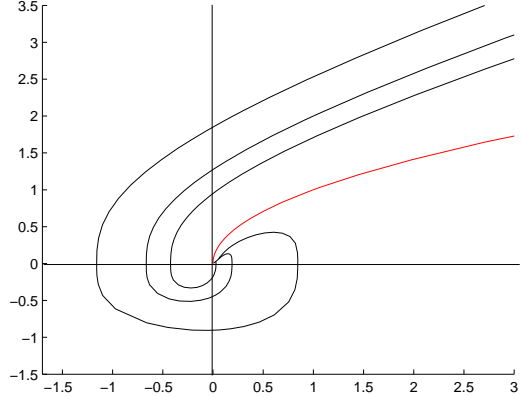}
   \end{center}
   \caption{Phase portrait around the origin for $c=1$. Experiment for $p=3$, $N=2$.}\label{figure1}
\end{figure}

\medskip

By remarks (a) and (c) above, it follows that for $c<1$, this unique
orbit disconnects from $Q_1$, hence it should cross at some point
the curve $U=cV+V^{p-1}$ (which still connects $P=(0,0)$ and $Q_1$);
as explained before, this orbit previously had a change of sign
(crossing the axis $U=0$) and then a change of monotonicity
(crossing the axis $V=0$). In particular, we can say that the
explicit orbit \eqref{expl.TW} is a separatrix between the monotone
and the non-monotone orbits. We draw the local phase portrait for
$c<1$, around the origin, in Figure~\ref{figure2} below. We gather
the discussion above in the following result.

\begin{lemma}\label{TWdim1}
(i) For any $c\in(0,1)$ and $K\ge 0$, there exists a unique
traveling wave solution
$\overline{f}_{c,K}(z)=\overline{f}_{c,K}(x-c\tau)$ of \eqref{a10}
in dimension $N=1$, having an interface at $z=K$ (that is,
$\overline{f}_{c,K}(z)=0$ for $z\ge K$) and moving with speed $c$.
In addition, $\overline{f}_{c,K}(z)=\overline{f}_{c,0}(z-K)$ for
$z\in\real$.

\noindent (ii) For $c=1$ and for any $K\ge 0$, there exists a unique
nonnegative traveling wave $f_{1,K}(z)=f_{1,K}(x-\tau)$ of
\eqref{a10} in dimension $N=1$ with interface at $z=K$, having the
explicit formula:
\begin{equation}\label{expl.TW2}
f_{1,K}(x-\tau)=\left(\frac{p-2}{p-1}\right)^{(p-1)/(p-2)}(K+\tau-x)_{+}^{(p-1)/(p-2)}.
\end{equation}
Here again, $f_{1,K}(z)=f_{1,0}(z-K)$ for $z\in\real$.

\noindent (iii) For any $c>1$ and $K\ge 0$, there exists a unique
traveling wave solution $f_{c,K}=f_{c,K}(x-c\tau)$ of \eqref{a10} in
dimension $N=1$ with interface at $z=K$ and moving with speed $c$.
Moreover, $f_{c,K}$ is nonnegative and non-increasing, and
$f_{c,K}(z)=f_{c,0}(z-K)$ for $z\in\real$.
\end{lemma}

\begin{figure}
   \begin{center}
   \includegraphics[width=8cm]{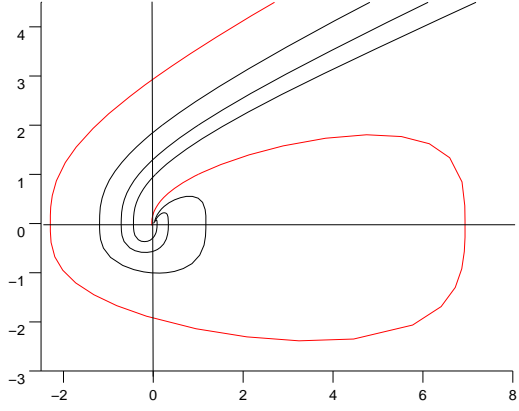}
   \end{center}
   \caption{Phase portrait around the origin for $c<1$. Experiment for $p=3$, $N=2$, $c=0.9$.}\label{figure2}
\end{figure}

\noindent \textbf{Compactly supported subsolutions for $0<c<1$.} We
are looking for nonnegative and compactly supported subsolutions
traveling with any speed $0<c<1$. These subsolutions are constructed
in the following way: from the analysis above, we know that, given
$c\in (0,1)$ and $K\ge 0$, there are two points $z_{c,K}\in
(-\infty,K)$ and $\tilde{z}_{c,K}\in (z_{c,K},K)$ such that
$$
z_{c,K} := \inf{\left\{ z\in (-\infty,K)\ :\ \overline{f}_{c,K}>0 \;\mbox{ in }\; (z,K) \right\}}>-\infty\,,
$$
and
$$
\overline{f}_{c,K}'>0 \;\mbox{ in }\; (z_{c,K},\tilde{z}_{c,K}) \;\;\;\mbox{ and }\;\;\; \overline{f}_{c,K}'<0 \;\mbox{ in }\; (\tilde{z}_{c,K},K)\,.
$$
We then define
\begin{equation}\label{subsol}
f_{c,K}(z)=\left\{\begin{array}{ll}\overline{f}_{c,K}(z), & \hbox{for} \ z_{c,K}\le z\le K, \\ 0, & \hbox{elsewhere}.\end{array}\right.
\end{equation}
In other words, we consider the positive hump of the graph of
$f_{c,K}$ located between its last change of sign and the interface.
It is immediate to check that $f_{c,K}$ is a compactly supported
subsolution to \eqref{a10} in dimension $N=1$, and that it has an
increasing part until reaching its maximum at $\tilde{z}_{c,K}$, and
then decreases  to the interface point $K$. The notation $f_{c,K}$
will designate in the sequel these subsolutions if $0<c<1$ and the
solutions to \eqref{a10} in dimension $N=1$ given by
Lemma~\ref{TWdim1} if $c\ge1$.

\subsection{Construction of subsolutions in dimension $N\ge 1$}\label{TWN}

We turn to equation \eqref{a10} posed in dimension $N\ge 1$ for
which we aim at constructing some special subsolutions having an
interface that moves out in all directions with a given velocity
$c<1$. The construction is based on the traveling waves $f_{c,K}$
identified in the previous subsection. The first attempt is to try
the form $V(\tau,x)=f_{c,K}(|x|-c\tau)$, $c\in (0,1)$, which
satisfies:
\begin{eqnarray*}
& & \partial_{\tau} V-\Delta_{p}V+|\nabla V|^{p-1}-V \\
& = &-cf_{c,K}'-\left(|f_{c,K}'|^{p-2}f_{c,K}'\right)'+|f_{c,K}'|^{p-1}-f_{c,K}-\frac{N-1}{|x|}|f_{c,K}'|^{p-2}f_{c,K}'\\
& \leq & -\frac{N-1}{|x|}|f_{c,K}'|^{p-2}f_{c,K}'\,.
\end{eqnarray*}
Thus, $V$ is a subsolution of \eqref{a10} in the region of $Q$ where
$f_{c,K}'\geq 0$. We therefore have to modify the profile in the
decreasing part of $f_{c,K}$ and we proceed as follows.

\medskip

\noindent \textbf{Traveling wave solutions to a modified equation in
dimension $N=1$}. For $\alpha\in (0,1/2)$, we consider the following
perturbation of \eqref{a10}:
\begin{equation}\label{a10m}
\partial_{\tau} \zeta -\partial_x \left( |\partial_x \zeta|^{p-2} \partial_x \zeta \right) + |\partial_x \zeta|^{p-1} -\a\,|\partial_x \zeta|^{p-2}\partial_x \zeta-\zeta=0\,, \quad (t,x)\in (0,\infty)\times\real\,,
\end{equation}
and look for traveling wave solutions $\zeta(\tau,x)=g(x-c\tau)$.
Then, $g$ solves
\begin{equation}\label{eq.modif}
-cg'-\left(|g'|^{p-2}g'\right)'+|g'|^{p-1}-\a\,|g'|^{p-2}g'-g=0.
\end{equation}
The phase-plane analysis for \eqref{eq.modif} is similar to that of
\eqref{OdeTW}, with the difference that an extra term
$-\a\,|V|^{p-2}V$ appears in the right-hand side of the second
equation in \eqref{syst1}. This is only reflected in the analysis at
infinity, where the point $Q_1$ changes into $(0,1/(1+\a))$ and the
explicit separatrix is obtained for $c=1/(1+\a)<1$. In particular,
we have the following analogue of Lemma~\ref{TWdim1}~(i).
\begin{lemma}\label{TWMdim1}
For any $\alpha>0$ sufficiently small, $c\in(0,1/(1+\a))$ and $K\ge
0$, there exists a unique traveling wave solution
$g_{c,K,\a}(z)=g_{c,K,\a}(x-c\tau)$ of \eqref{a10m} having an
interface at $z=K$ and moving with speed $c$. In addition,
$g_{c,K,\a}(z)=g_{c,0,\a}(z-K)$ for $z\in\real$ and there are two
points $z_{c,K,\a}\in (-\infty,K)$ and $\tilde{z}_{c,K,\a}\in
(z_{c,K,\a},K)$ such that
$$
z_{c,K,\a} := \inf{\left\{ z\in (-\infty,K)\ :\ g_{c,K,\a}>0 \;\mbox{ in }\; (z,K) \right\}}>-\infty\,,
$$
and
$$
g_{c,K,\a}'>0 \;\mbox{ in }\; (z_{c,K,\a},\tilde{z}_{c,K,\a}) \;\;\;\mbox{ and }\;\;\; g_{c,K,\a}'<0 \;\mbox{ in }\; (\tilde{z}_{c,K,\a},K)\,.
$$
\end{lemma}
Setting
$$
M_{c,\a} :=  \sup_{z\in [z_{c,0,\a},0]}{\{ g_{c,0,\a}(z)\}}\,,
$$
we notice that
\begin{equation}\label{gaston}
z_{c,K,\a}=z_{c,0,\a}+K\,, \quad
\tilde{z}_{c,K,\a}=\tilde{z}_{c,0,\a}+K\,, \quad \sup_{z\in
[z_{c,K,\a},K]}{\{ g_{c,K,\a}(z)\}} = M_{c,\alpha}\,.
\end{equation}
If we put now $V(\tau,x)=g_{c,K,\a}(|x|-c\tau)$, we calculate and find that
\begin{equation*}
\partial_{\tau} V-\Delta_{p}V+|\nabla V|^{p-1}-V=\left(\a-\frac{N-1}{|x|}\right)\ \left( |g_{c,K,\a}'|^{p-2}g_{c,K,\a}' \right)(|x|-c\tau)\,,
\end{equation*}
and it is a subsolution where $g_{c,K,\a}'\leq 0$ and
$\a\geq(N-1)/|x|$. Matching these two conditions turns out to be
possible as we show now.

Fix $c\in (1/2,1)$ and $\a_c:=(1-c)/(1+c)$ and define
\begin{equation}\label{jeanne}
\tau_0(c) := \max{\left\{ \frac{2(N-1)}{\alpha_c} - 2 \tilde{z}_{c,0,\a_c} , - \frac{\tilde{z}_{c,0,\a_c}}{c}  \right\}} > \frac{2(N-1)}{\alpha_c}\,,
\end{equation}
the point $\tilde{z}_{c,0,\a_c}\in (-\infty,0)$ being defined in
Lemma~\ref{TWMdim1}. Then $c<1/(1+\a_c)$ and, for $K\ge 0$, $\tau\ge
\tau_0(c)$, and $|x|\ge \tilde{z}_{c,K,\a_c} +
c\tau=\tilde{z}_{c,0,\a_c}+K+c\tau$, we have
$$
\frac{N-1}{|x|} \le \frac{N-1}{\tilde{z}_{c,0,\a_c}+c\tau_0(c)} \le \frac{2(N-1)}{2\tilde{z}_{c,0,\a_c}+\tau_0(c)} \le \alpha_c\,,
$$
and
\begin{eqnarray*}
g_{c,K,\a_c}'(|x|-c\tau) < 0 & \mbox{ if } & \tilde{z}_{c,K,\a_c} + c\tau \le |x| < K+c\tau\,, \\
g_{c,K,\a_c}'(|x|-c\tau) = 0 & \mbox{ if } & K + c\tau \le |x|\,.
\end{eqnarray*}
Consequently, for $c\in (1/2,1)$, $\alpha_c=(1-c)/(1+c)$, and $K>0$,
the function $V$ defined by $V(\tau,x)=g_{c,K,\a_c}(|x|-c\tau)$ is a
subsolution to \eqref{a10} for $\tau\ge \tau_0(c)$ and $|x|\ge
\tilde{z}_{c,K,\a_c} + c\tau$. Observing that any positive constant
is a subsolution to \eqref{a10}, we construct a compactly supported
subsolution $v_{c,K}$ to \eqref{a10} by setting
\begin{equation}\label{subs.dimN}
v_{c,K}(\tau,x):=\left\{\begin{array}{ll}
M_{c,\alpha_c} & \quad \hbox{if} \ 0\leq|x|<\tilde{z}_{c,K,\a_c} + c\tau\,,\\
& \\
g_{c,K,\alpha_c}(|x|-c\tau) & \quad \hbox{if} \ |x|>\tilde{z}_{c,K,\a_c} + c\tau\,,
\end{array}\right.
\end{equation}
for $\tau\ge \tau_0(c)$. It is easy to check that the function
$v_{c,K}$ is a subsolution to \eqref{a10} in
$[\tau_0(c),\infty)\times\real^N$. It will be used for comparison
from below, as indicated in the next subsection.

\subsection{Proof of Theorem~\ref{main.posit}}\label{subsect.mainpos}

We conclude the proof of Theorem~\ref{main.posit} by a comparison
argument, using the subsolutions and supersolutions constructed in
the previous subsections. Before that, we identify a class of
solutions of \eqref{a10} that is representative for the general
solutions.

We say that a function $V=V(\tau,x)$ is \emph{radially
non-increasing} if $V(\tau,\cdot)$ is radially symmetric for all
$\tau$, and it is non-increasing in the radial variable $r:=|x|$.
For example, the subsolutions $v_{c,K}$ are radially non-increasing.
The next results show that the class of radially non-increasing
solutions of \eqref{a10} is sufficient for our aims.
\begin{lemma}\label{rad.noninc}
Let $u_0=u_0(r)$ be a radially non-increasing function satisfying
\eqref{a3}. Then, the solution $v$ of \eqref{a10} with initial
condition $u_0$ is also radially non-increasing.
\end{lemma}

\noindent\textbf{Proof.} The radial symmetry of the solution $v$
follows from the invariance of the equation \eqref{a10} with respect
to rotations. We write now the equation satisfied by $\xi=\partial_r
v$, obtained by differentiating \eqref{a10} with respect to $r$:
$$
\partial_t\xi-\partial_r^2(|\xi|^{p-2}\xi)-\frac{N-1}{r}\partial_r(|\xi|^{p-2}\xi)+\frac{N-1}{r^2}|\xi|^{p-2}\xi+(p-1)|\xi|^{p-3}\xi \partial_r\xi-\xi=0,
$$
which is a parabolic equation (of porous medium type) and satisfies
a maximum principle. Since $0$ is a solution to the above equation,
the derivative $\xi=\partial_r v$ remains nonpositive if it is
initially nonpositive and it follows that $v$ is radially
non-increasing. \qed

\medskip

We are now in position to end the proof of Theorem~\ref{main.posit}
for radially non-increasing initial data. More precisely, we have
the following upper and lower bounds for the edge $\varrho\,(\tau)$
defined in \eqref{bs32} of the support of $v(\tau)$.
\begin{lemma}\label{leb3}
Let $u_0=u_0(r)$ be a radially non-increasing function satisfying
\eqref{a3} and denote the solution of \eqref{a10} with initial
condition $u_0$ by $v$. For any $c\in (1/2,1)$, there exists
$\tau_1(c)>0$ such that, for any $\tau\geq\tau_1(c)$, we have:
\begin{equation}\label{b4}
1+c(\tau-\tau_1(c)) \leq\varrho\,(\tau) \le R_0 + \frac{p-1}{p-2}\ \Vert u_0\Vert_\infty^{(p-2)/(p-1)} + \tau\,.
\end{equation}
In particular, we obtain that $\varrho(\tau)/\tau\to1$ as
$\tau\to\infty$.
\end{lemma}

\noindent\textbf{Proof.} The upper bound follows by comparison with
the explicit traveling wave solutions \eqref{expl.TW2}. More
precisely, we define
\begin{equation}\label{bs41}
R_1:= R_0 + \frac{p-1}{p-2}\ \Vert u_0\Vert_\infty^{(p-2)/(p-1)}
\end{equation}
and consider the function
$\overline{v}(\tau,x)=f_{1,R_1}(x_1-\tau)$, which is a solution of
\eqref{a10} by Lemma~\ref{TWdim1}. If $x=(x_i)_{1\le i \le N}\in\RR^N$ is such that
$x_1\ge R_0$, then $|x|\ge R_0$ and $u_0(x)=0\le\overline{v}(0,x)$
while, if $x_1\le R_0$, \bean
u_0(x) & \le & \Vert u_0\Vert_\infty \le \left( \frac{p-2}{p-1} \right)^{(p-1)/(p-2)}\ \left( R_1 - R_0 \right)^{(p-1)/(p-2)} \\
& \le & \left( \frac{p-2}{p-1} \right)^{(p-1)/(p-2)}\ \left( R_1 -
x_1 \right)^{(p-1)/(p-2)} = \overline{v}(0,x)\,. \eean The
comparison principle then entails that $v(\tau,x)\le
\overline{v}(\tau,x)$ for $(\tau,x)\in [0,\infty)\times\real^N$,
from which we conclude that $\mathcal{P}_v(\tau) \subseteq \left\{
x\in\RR^N\ : \ x_1\le R_1+\tau \right\}$. Owing to the rotational
invariance of \eqref{a10}, we actually have $\mathcal{P}_v(\tau)
\subseteq \left\{ x\in\RR^N\ : \ \langle x , \omega \rangle \le
R_1+\tau \right\}$ for every $\omega\in\mathbb{S}^{N-1}$ and
$\tau\ge 0$, and thus
\begin{equation}\label{fantasio}
\mathcal{P}_v(\tau)\subseteq B(0,R_1+\tau)\,.
\end{equation}

The lower bound follows from comparison with the subsolutions
constructed in \eqref{subs.dimN}.  Fix $c\in (1/2,1)$
and put $r_1:=1+c\tau_0(c)$, $\tau_0(c)$ being defined by
\eqref{jeanne}. Since $v(\tau)$ is radially non-increasing for all
$\tau\ge 0$ by Lemma~\ref{rad.noninc}, we infer from
Proposition~\ref{prb1} that, for $x\in B(0,r_1)$ and $\tau\ge
T_{r_1}$,
$$
v(\tau,x) \ge v\left( \tau, \frac{r_1 x}{|x|} \right) \ge \varepsilon_{r_1}\ \tau^{(p-1)/(p-2)}\,.
$$
Define $\tau_1(c)$ by
$$
\tau_1(c) := \max{\left\{ \tau_0(c), T_{r_1} , \left( \frac{M_{c,(1-c)/(1+c)}}{\varepsilon_{r_1}} \right)^{(p-2)/(p-1)}\right\}}\,,
$$
so that the previous inequality and the properties of $v_{c,1}$
defined in \eqref{subs.dimN} guarantee that
$$
v(\tau_1(c),x) \ge M_{c,(1-c)/(1+c)} \ge v_{c,1}(\tau_0(c),x)\,, \quad x\in B(0,r_1)\,.
$$
Since $v_{c,1}(\tau_0(c),x)=0$ for $x\not\in B(0,r_1)$, we also have
$v(\tau_1(c),x) \ge v_{c,1}(\tau_0(c),x)$ for $x\not\in B(0,r_1)$.
Recalling that $v_{c,1}$ is a subsolution to \eqref{a10} in
$(\tau_0(c),\infty)\times\real^N$, we infer from the comparison
principle that
\begin{equation}\label{compTW}
v(\tau+\tau_1(c),x)\geq v_{c,1}(\tau+\tau_0(c),x), \quad
(\tau,x)\in Q\,.
\end{equation}
Consequently, $v(\tau+\tau_1(c),x)>0$ if $x\in B(0,r_1+c\tau)$, whence
\begin{equation}\label{lebrac}
B(0,1+c(\tau+\tau_0(c)-\tau_1(c))) \subset \mathcal{P}_v(t)\,, \quad \tau\ge \tau_1(c)\,.
\end{equation}
This readily implies that
$$
\varrho(\tau) \ge 1+c(\tau+\tau_0(c)-\tau_1(c)) \ge 1+c(\tau-\tau_1(c))\,, \quad \tau\ge \tau_1(c)\,.
$$
In particular, we deduce from \eqref{fantasio} and \eqref{lebrac}
that
\begin{equation*}
\liminf\limits_{\tau\to\infty}\frac{\varrho\,(\tau)}{\tau}\ge c \;\;\;\mbox{ for any }\;\;\; c\in (1/2,1) \;\;\;\mbox{ and }\;\;\;
\limsup\limits_{\tau\to\infty}\frac{\varrho\,(\tau)}{\tau}\le 1\,,
\end{equation*}
which implies that $\varrho(\tau)/\tau\to1$ as $\tau\to\infty$.
\qed

\medskip
Rephrasing the rescaling and coming back to the notation with
$t=(e^{(p-2)\tau}-1)/(p-2)$ and $\gamma(t)=\varrho(\tau)$, we find
the result of Theorem~\ref{main.posit} for radially non-increasing
inital data. The extension to arbitrary initial data satisfying
\eqref{a3} is performed in Section~\ref{sec:aid}. Moreover, we
notice that the speed is the same in any direction
$\omega\in\mathbb{S}^{N-1}$, as stated.


\section{Proof of Theorem~\ref{asympt.main}}\label{Sect.scaling2}

\subsection{Scaling variables~II}\label{subsec4.1}

According to Proposition~\ref{prb1}, as $\tau\to\infty$ the solution
$v$ to \eqref{a10}, \eqref{a11} expands in space and grows
unboundedly in time. In order to take into account such phenomena,
we introduce next a further scaling of the space variable
\begin{equation}\label{def.y}
y:= \frac{x}{1+\tau} \,,
\end{equation}
together with the new unknown function $w=w(\tau,y)$ defined by
\begin{equation}\label{def.w}
v(\tau,x) =(1+\tau)^{(p-1)/(p-2)}\ w\left(\tau,\frac{x}{1+\tau} \right)\,, \quad (\tau,x)\in [0,\infty)\times\real^N\,.
\end{equation}
It follows from \eqref{a10} and \eqref{a11} that $w$ solves
\begin{equation}\label{c1}
\partial_{\tau}w -\frac{1}{1+\tau}\left(\Delta_p w + y \cdot \nabla w -\frac{p-1}{p-2}\,w\right)+|\nabla w|^{p-1} -w  = 0\ \,,
\quad (\tau,y)\in Q\,,
\end{equation}
with the same initial condition
\begin{equation}\label{c2}
w(0) = u_0\,, \quad y\in\RR^N\,.
\end{equation}
Throughout this section we assume that $u_0$ is radially
non-increasing besides \eqref{a3}. In particular, $u_0(0)>0$. We
gather several properties of $w$ in the next lemma.
\begin{lemma}\label{lec1}
There is a positive constant $C_4$ depending only on $p$, $N$, and $u_0$ such that
\begin{equation} \label{c3}
\Vert w(\tau)\Vert_1 + \Vert w(\tau)\Vert_\infty + \Vert \nabla w(\tau)\Vert_\infty \le C_4 \,, \quad \tau\ge 0\,,
\end{equation}
\begin{equation} \label{c3b}
w(\tau,y) \ge \frac{1}{C_4} \left( r_*^2 - |y|^2 \right)_+^{(p-1)/(p-2)}\,, \quad (\tau,y)\in Q\,,
\end{equation}
the radius $r_*$ being defined in Corollary~\ref{cor:bfbsb}.
Moreover,
\begin{equation}\label{c4}
\mathcal{P}_w(\tau):= \left\{ y\in\RR^N : \;\; w(\tau,y)>0 \right\} \subseteq B\left(0, 1 + \frac{R_1}{1+\tau}\right)
\end{equation}
for $\tau\ge 0$ where $R_1$ is defined by \eqref{bs41}. In addition,
for any $c\in (1/2,1)$, we have
\begin{equation}\label{c4b}
B\left( 0, c - \frac{\tau_1(c)}{1+\tau} \right)\subset \mathcal{P}_{w}(\tau) \quad {\rm for } \ \tau\geq \tau_1(c)\,,
\end{equation}
the time $\tau_1(c)>0$ being defined in Lemma~\ref{leb3}.
\end{lemma}

\noindent\textbf{Proof.} The estimates \eqref{c3} and \eqref{c3b}
readily follow from \eqref{a12} and \eqref{spirou}, while \eqref{c4}
is a consequence of \eqref{fantasio}. The assertion about the ball
$B(0,c-\tau_1(c)/(1+\tau))$ follows from \eqref{lebrac}.\qed

\medskip

At this point, \eqref{c1} indicates that $w(\tau)$ behaves as
$\tau\to\infty$ as the solution $\tilde{w}$ to the Hamilton-Jacobi
equation $\partial_\tau \tilde{w} + |\nabla\tilde{w}|^{p-1} -
\tilde{w}=0$ in $Q$ which is known to converge to a stationary
solution uniquely determined by the limit of the support of
$\tilde{w}(\tau)$ as $\tau\to\infty$, see, e.g.,
\cite[Theorem~A.2]{La08}. As an intermediate step, we thus have to
identify the limit of the support of $w(\tau)$ as $\tau\to\infty$.
Thanks to \eqref{c4}, we already know that it is included in
$B(0,1)$ but the information in \eqref{c4b}  are yet too weak to
exclude the vanishing of $w(\tau)$ outside a ball of radius smaller
than one. To complete the proof of Theorem~\ref{asympt.main} for
radially non-increasing initial data, we show first that the
asymptotic limit is supported exactly in the ball $B(0,1)$. Then we
use a viscosity technique, the same that has been used in the
previous paper \cite{LV07} to establish the convergence to the
expected stationary solution.

\subsection{Proof of Theorem~\ref{asympt.main}: $N=1$}\label{subsec4.2}

We first consider the one-dimensional case $N=1$ and divide the proof into several technical steps.

\medskip

\noindent \textbf{Step 1. A special family of subsolutions}. Given
$c\in (1/2,1)$, we have
$$
v(\tau,x)\ge v_{c,1}(\tau+\tau_0(c)-\tau_1(c),x)\,, \quad (\tau,x)\in [\tau_1(c),\infty)\times\real\,,
$$
by \eqref{compTW}, the times $\tau_0(c)$ and $\tau_1(c)$ being defined in \eqref{jeanne} and Lemma~\ref{leb3}, respectively. Then,
\begin{equation}\label{compTW2}
w(\tau,y) \ge w_c(\tau,y) := \frac{1}{(1+\tau)^{(p-1)/(p-2)}}\ v_{c,1}(\tau+\tau_0(c)-\tau_1(c),y(1+\tau))
\end{equation}
for $(\tau,y)\in [\tau_1(c),\infty)\times\real$.

\medskip

\noindent \textbf{Step 2. An explicit family of supersolutions}. Let us introduce the following family of functions:
\begin{equation}\label{FR1}
F_{R}(\tau,y)=\left(\frac{p-2}{p-1}\right)^{(p-1)/p-2)} \left(\frac{\tau+R}{\tau+1}-|y|\right)_{+}^{(p-1)/(p-2)}\,, \quad (\tau,y)\in Q\,.
\end{equation}
We easily obtain by direct calculation that $F_R$ is a classical
solution of \eqref{c1} for $y\neq 0$, and for all parameter values
$R\ge 0$. However, near $y=0$, it is only a supersolution both in
the weak and the viscosity sense. The latter is straightforward to
verify using the definition of viscosity subsolutions and
supersolutions with jets, as in the classical survey \cite{CIL92}.
Let us mention at this point that these functions can be used in a
comparison argument to give an alternative proof of \eqref{c4}.

\begin{remark}\label{rem:nat} This family of functions arises naturally if we think about asymptotics.
Indeed, as already mentioned, we formally expect that the asymptotic
profiles of \eqref{c1} should be given by solutions of the
stationary Hamilton-Jacobi equation
\begin{equation}\label{limit}
|\nabla \tilde{w}|^{p-1}-\tilde{w}=0,
\end{equation}
supported in some ball $B(0,R)$, that is
\begin{equation*}
H_{R}(y):=\left(\frac{p-2}{p-1}\right)^{(p-1)/(p-2)}\left(R-|y|\right)_{+}^{(p-1)/(p-2)}\,, \quad y\in\real\,.
\end{equation*}
Making then the ``ansatz" that, for large times, the solution of \eqref{c1} should behave in a similar way as its limit, we write
\begin{equation*}
w(\tau,y)\sim \left(\frac{p-2}{p-1}\right)^{(p-1)/(p-2)}\left(C(\tau)-|y|\right)_{+}^{(p-1)/(p-2)}\,.
\end{equation*}
Integrating the resulting ordinary differential equation for
$C(\tau)$, we arrive at the family of explicit exact profiles $F_R$
given by \eqref{FR1}.
\end{remark}

\medskip

\noindent \textbf{Step 3. Constructing suitable subsolutions.} We
now face the problem of finding suitable subsolutions with similar
behaviour. Since the $F_R$'s are classical solutions to \eqref{c1}
except at $y=0$, we expect to be able to construct also a family of
subsolutions based on them. To this end, we consider the ``damped''
family $F_{R,\vartheta,\beta}$ defined by
\begin{equation}\label{marsupilami}
F_{R,\vartheta,\beta}(\tau,y) := \vartheta\ \left( \frac{p-2}{p-1} \right)^{(p-1)/(p-2)} \left(  \frac{\beta(\tau+R)}{\tau+1} - |y| \right)_+^{(p-1)/(p-2)}\,, \quad (\tau,y)\in Q\,,
\end{equation}
for parameters $R\in (0,1)$, $\vartheta\in (0,1]$, and $\beta\in (1/2,1]$. Observe that, since $(p-1)/(p-2)>1$, $F_{R,\vartheta,\beta}$ and $|\nabla F_{R,\vartheta,\beta}|^{p-2} \nabla F_{R,\vartheta,\beta}$ both belong to $\mathcal{C}^1([0,\infty)\times(\real\setminus\{0\}))$. For $\vartheta\in (0,1)$, $\beta\in (1/2,1]$, $\tau>0$ and $y\neq0$, we calculate
\begin{equation*}
\begin{split}
\partial_\tau F_{R,\vartheta,\beta}&-\frac{1}{1+\tau}\left(\Delta_{p}F_{R,\vartheta,\beta}+y\cdot\nabla F_{R,\vartheta,\beta}-\frac{p-1}{p-2}F_{R,\vartheta,\beta}\right)+|\nabla F_{R,\vartheta,\beta}|^{p-1}-F_{R,\vartheta,\beta}\\
&=\vartheta \beta \frac{1-R}{(1+\tau)^2} F_{R,1,\beta}^{1/(p-1)} - \frac{\vartheta}{1+\tau} \left( \vartheta^{p-2} - \frac{\beta (\tau+R)}{\tau+1} \right)\ F_{R,1,\beta}^{1/(p-1)} - \vartheta (1-\vartheta^{p-2}) F_{R,1,\beta}\\
&=\vartheta \left( \frac{\beta-\vartheta^{p-2}}{1+\tau} - (1-\vartheta^{p-2}) F_{R,1,\beta}^{(p-2)/(p-1)} \right)\ F_{R,1,\beta}^{1/(p-1)}\\
&\le \vartheta(1-\vartheta^{p-2}) F_{R,1,\beta}^{1/(p-1)} \left[\frac{1}{1+\tau}- \frac{p-2}{p-1}\left(\frac{\beta(\tau+R)}{\tau+1}-|y|\right)\right].
\end{split}
\end{equation*}
Analyzing the sign of the last expression and taking into account that $\vartheta\in(0,1)$, we obtain that $F_{R,\vartheta,\beta}$ has the following properties:
\begin{equation}\label{pim}
\begin{minipage}{10cm}
$F_{R,\vartheta,\beta}$ is a classical subsolution to \eqref{c1} in \\ $\{ (\tau,y)\in Q \ : \ \tau\ge\tau_2(R,\beta)\,, \ 0<|y|\le K_{R,\beta}(\tau) \}$
\end{minipage}
\end{equation}
with
\begin{equation}\label{pam}
\tau_2(R,\beta) := \frac{p-1}{\beta (p-2)} - R \;\;\mbox{ and }\;\;  K_{R,\beta}(\tau) := \frac{\beta(\tau+R)}{\tau+1} - \frac{p-1}{p-2}\ \frac{1}{\tau+1}\,,
\end{equation}
and
\begin{equation}\label{poum}
F_{R,\vartheta,\beta} \;\mbox{ vanishes for }\; |y|\geq\frac{\beta(\tau+R)}{\tau+1} \;\mbox{ and }\; \tau\ge 0\,.
\end{equation}
Let us notice here that both the edge of the support of $F_{R,\vartheta,\beta}$ and the constant $K_{R,\beta}(\tau)$, where the behaviour changes, do not depend on $\vartheta$. While the two properties \eqref{pim} and \eqref{poum} are suitable for our purpose, the function $F_{R,\vartheta,\beta}$ does not behave in a suitable way near $y=0$ (where it is a viscosity supersolution) and in an asymptotically small region near the edge of its support (where it is a classical supersolution). However, we already have a positive bound from below for $w$ in a small neighbourhood of $y=0$ by \eqref{c3b} which allows us to remedy to the first bad property of $F_{R,\vartheta,\beta}$. More precisely, we infer from \eqref{c3b} that
$$
w(\tau,y) \ge C_5 := \frac{1}{C_4} \left( \frac{3r_*^2}{4} \right)^{(p-1)/(p-2)}>0\,, \quad (\tau,y)\in [0,\infty)\times B(0,r_*/2)\,,
$$
whence
\begin{equation}\label{ld15}
w(\tau,y) \ge \vartheta \ge F_{R,\vartheta,\beta}(\tau,y) \,, \quad (\tau,y)\in [0,\infty)\times B(0,r_*/2)\,,
\end{equation}
provided that
\begin{equation}\label{ld16}
0<\vartheta < \min{\{1,C_5\}}\,.
\end{equation}
Consider next
$$
\tau\ge\tau_2(R,\beta) \;\;\mbox{ and }\;\; K_{R,\beta}(\tau)\le |y| \le \frac{\beta(\tau+R)}{\tau+1}\,.
$$
Then
\begin{eqnarray}
F_{R,\vartheta,\beta}(\tau,y) & \le & \vartheta\ \left( \frac{p-2}{p-1} \right)^{(p-1)/(p-2)}\ \left( \frac{p-1}{p-2}\ \frac{1}{1+\tau} \right)^{(p-1)/(p-2)} \nonumber \\
& = & \frac{\vartheta}{(1+\tau)^{(p-1)/(p-2)}}\,.\label{ld17}
\end{eqnarray}
Now, if $c\in (\beta,1)$, we have
$$
|y| (1+\tau) \le \beta (\tau+R) \le \tilde{z}_{c,1,(1-c)/(1+c)} + c (\tau+\tau_0(c)-\tau_1(c))
$$
as soon as
\begin{equation}\label{ld18}
\tau\ge \tau_3(c,R,\beta) := \frac{\beta R+c(\tau_1(c)-\tau_0(c))-\tilde{z}_{c,1,(1-c)/(1+c)}}{c-\beta}\,.
\end{equation}
In that case,
$$
w_c(\tau,y) = \frac{1}{(1+\tau)^{(p-1)/(p-2)}}\ v_{c,1}(\tau+\tau_0(c)-\tau_1(c),y(1+\tau)) = \frac{M_{c,(1-c)/(1+c)}}{(1+\tau)^{(p-1)/(p-2)}}
$$
according to the properties \eqref{subs.dimN} of $v_{c,1}$.
Recalling \eqref{compTW2} and \eqref{ld17} we realize that
\begin{equation}\label{ld19}
F_{R,\vartheta,\beta}(\tau,y) \le w_c(\tau,y) \le w(\tau,y)\,, \quad K_{R,\beta}(\tau)\le |y| \le \frac{\beta(\tau+R)}{\tau+1}\,,
\end{equation}
provided
\begin{equation}\label{ld20}
c\in (\beta,1)\,, \quad \vartheta<\min{\{1, M_{c,(1-c)/(1+c)}\}}\,, \quad \tau\ge \max{\{ \tau_1(c) , \tau_2(R,\beta), \tau_3(c,R,\beta)\}}\,.
\end{equation}

After this preparation, we are in a position to establish a positive
lower bound for $w$ on the ball $B(0,1-\e)$ for any $\e\in (0,1/4)$.
Indeed, we fix $\e\in (0,1/4)$, choose $c=1-\e$, $R=\beta=1-2\e$,
and define
$$
\tau_4(\e) := \max{\left\{ \frac{\tau_1(1-\e)}{\e} , \tau_2(1-2\e,1-2\e), \tau_3(1-\e,1-2\e,1-2\e) \right\}}\,.
$$
As $\tau_4(\e)>\tau_1(1-\e)/\e$, \eqref{c4b} guarantees that
$B(0,1-2\e)\subset \mathcal{P}_w(\tau_4(\e))$ and there is thus
$m_\e\in (0,1)$ such that
\begin{equation}\label{ld21}
w(\tau_4(\e),y)\ge m_\e\,, \quad y\in B(0,1-2\e)\,.
\end{equation}
Now, for $\vartheta\in (0,1)$ satisfying
\begin{equation}\label{ld22}
0<\vartheta<\min{\{ m_\e, C_5, M_{1-\e,\e/(2-\e)} \}}
\end{equation}
we infer from \eqref{pam}, \eqref{ld15}, \eqref{ld16}, \eqref{ld19}, \eqref{ld20}, and \eqref{ld21} that
$$
F_{1-2\e,\vartheta,1-2\e}(\tau,y) \le w(\tau,y)\,, \quad |y|\in\left\{ \frac{r_*}{2} , K_{1-2\e,1-2\e}(\tau) \right\}\,, \quad \tau\ge \tau_4(\e)\,,
$$
and
$$
F_{1-2\e,\vartheta,1-2\e}(\tau_4(\e),y) \le \vartheta \le m_\e \le w(\tau_4(\e),y)\,, \quad \frac{r_*}{2} \le |y| \le K_{1-2\e,1-2\e}(\tau_4(\e)) \le 1-2\e\,.
$$
It then follows from \eqref{c1}, \eqref{pim}, and the comparison principle that
$$
F_{1-2\e,\vartheta,1-2\e}(\tau,y) \le w(\tau,y)\,, \quad \frac{r_*}{2} \le |y| \le K_{1-2\e,1-2\e}(\tau) \,, \quad \tau\ge \tau_4(\e)\,.
$$
Recalling \eqref{poum}, \eqref{ld15}, and \eqref{ld19}, we have thus established that
\begin{equation}\label{ld23}
F_{1-2\e,\vartheta,1-2\e}(\tau,y) \le w(\tau,y)\,, \quad \tau\in [\tau_4(\e),\infty)\times\real\,,
\end{equation}
for all $\vartheta\in (0,1)$ satisfying \eqref{ld22}.

\medskip

\noindent \textbf{Step 4. Positive bound from below.} For $\e\in
(0,1/4)$, fix $\vartheta_\e\in (0,1)$ satisfying \eqref{ld22}.
According to \eqref{ld23}, we have, for $\tau\ge\tau_4(\e)+1$ and
$y\in B(0,1-3\e)$,
\begin{eqnarray*}
w(\tau,y) & \ge & \vartheta_\e \left( \frac{p-2}{p-1} \right)^{(p-1)/(p-2)}\ \left( \frac{(1-2\e) (\tau+1-2\e)}{\tau+1} - |y| \right)_+^{(p-1)/(p-2)} \\
& \ge & \vartheta_\e \left( \frac{p-2}{p-1} \right)^{(p-1)/(p-2)}\ \left( \frac{\e (\tau-1+4\e)}{\tau+1} \right)_+^{(p-1)/(p-2)}\\
& \ge & \mu_\e:=\vartheta_\e \left( \frac{2(p-2) \e^2}{p-1} \right)^{(p-1)/(p-2)}>0\,.
\end{eqnarray*}
We have thus proved that, for all $\e\in (0,1/4)$, there are $\mu_\e>0$ and $\tau_5(\e):=\tau_4(\e)+1$ such that
\begin{equation}\label{ld24}
0 < \mu_\e \le w(\tau,y)\,, \quad (\tau,y)\in [\tau_5(\e),\infty)\times B(0,1-3\e)\,.
\end{equation}

\medskip

\noindent \textbf{Step 5. Convergence. Viscosity argument.} To
complete the proof, we use an argument relying on the theory of
viscosity solutions in a similar way as in the paper \cite{LV07} for
the subcritical case of \eqref{a1} with $q\in(1,p-1)$. We thus
employ the technique of half-relaxed limits \cite{BlP88} in the same
fashion as in \cite[Section~3]{Ro01} and \cite{LV07}. To this end,
we pass to the logarithmic time and introduce the new variable
$s:=\log(1+\tau)$ along with the new unknown function
$$
w(\tau,y) = \omega(\log{(1+\tau)} , y)\,, \quad (\tau,y)\in [0,\infty)\times\real\,.
$$
Then, $\partial_\tau w(\tau,y)=e^{-s} \partial_s\omega(s,y)$ and it
follows from \eqref{c1} and \eqref{c2} that $\omega$ solves
\begin{equation}\label{log-time}
e^{-s} \left( \partial_s \omega - \Delta_p \omega - y \cdot \nabla \omega + \frac{p-1}{p-2}\ \omega \right) + \vert\nabla \omega \vert^{p-1} - \omega  = 0\ \,, \quad (s,y)\in Q\,,
\end{equation}
with initial condition $\omega(0)=u_0$. We readily infer from Lemma~\ref{lec1} that
\begin{eqnarray}
& & \Vert \omega(s)\Vert_1 + \Vert \omega(s)\Vert_\infty + \Vert \nabla \omega(s)\Vert_\infty \le C_4 \,, \quad s\ge 0\,, \label{cvld20} \\
& & \omega(s,y) = 0 \;\;\;\mbox{ for }\;\;\; s\ge 0 \;\;\;\mbox{ and }\;\;\; |y|\ge 1+R_1 e^{-s}\,. \label{cvld21}
\end{eqnarray}
We next introduce the half-relaxed limits
$$
\omega_*(y) := \liminf_{(\sigma,z,\lambda) \to (\s,y,\infty)}{\omega(\lambda+\sigma,z)} \quad\mbox{ and }\quad \omega^*(y) := \limsup_{(\sigma,z,\lambda) \to (\s,y,\infty)}{\omega(\lambda+\sigma,z)},
$$
for $(s,y)\in Q$, which are well-defined according to the uniform
bounds in \eqref{cvld20} and indeed do not depend on $s>0$. Then,
the definition of $\omega_*$ and $\omega^*$ clearly ensures that
\begin{equation} \label{c6}
0 \le \omega_*(y) \le \omega^*(y) \;\;\mbox{ for }\;\; y\in\RR \,,
\end{equation}
while the uniform bounds \eqref{cvld20} and the Rademacher theorem
warrant that $\omega_*$ and $\omega^*$ both belong to
$W^{1,\infty}(\RR)$. Finally, by Proposition~\ref{pert} applied to
\eqref{log-time}, $\omega_*$ and $\omega^*$ are viscosity
supersolution and subsolution, respectively, to the Hamilton-Jacobi
equation
\begin{equation} \label{c7}
H(\zeta,\nabla\zeta) :=|\nabla\zeta|^{p-1} - \zeta = 0 \;\;\;\mbox{ in }\;\;\; \real\,.
\end{equation}

Our aim is now to show that $\omega_*\ge \omega^{*}$ in $\real$
(which implies that $\omega_*=\omega^{*}$ by \eqref{c6}). Since
$\omega^{*}$ and $\omega_*$ are subsolution and supersolution to
\eqref{c7}, respectively, such an inequality would follow from a
comparison principle which cannot be applied yet without further
information on $\omega^{*}$ and $\omega_*$. We actually need to
prove the following two facts:
\begin{itemize}
\item[(a)] $\omega_*(y)=\omega^*(y)=0$ if $|y|\ge 1$,
\item[(b)] $\omega^{*}(y)\geq \omega_*(y)>0$ if $y\in B(0,1)$,
\end{itemize}
and then to follow the technique used in \cite{LV07} to conclude
that $\omega_*=\omega^*$ and identify the limit.

To prove assertion~(a), let us take $y\in\real$ with $|y|>1$. We
then deduce from \eqref{cvld21} that there exists $s_1(y)>0$ such
that $\omega(s,y)=0$ for $s\ge s_1(y)$. Pick sequences
$(\sigma_n)_{n\geq1}$, $(\lambda_n)_{n\geq1}$, and $(z_n)_{n\geq1}$
such that $\sigma_n\to0$, $\lambda_n\to\infty$, $z_n\to y$, and
$\omega(\sigma_n+\lambda_n,z_n)\to \omega^{*}(y)$. On the one hand,
there exists $n_1(y)>0$ such that $\sigma_n+\lambda_n>s_1(y)$ for
any $n\geq n_1(y)$; hence $\omega(\sigma_n + \lambda_n,y) = 0$ for
any $n\geq n_1(y)$. On the other hand, we can write:
$$
|\omega(\sigma_n+\lambda_n,z_n)-\omega(\sigma_n+\lambda_n,y)|\leq|y-z_n|\|\nabla \omega(\sigma_n+\lambda_n)\|_{\infty}\leq C_{4}|y-z_n|\to0,
$$
hence $\omega^{*}(y)=0=\omega_{*}(y)$ for any $y\in\real$ with
$|y|>1$. In addition, since $\omega^{*}$ and $\omega_{*}$ are
continuous, it follows that $\omega^{*}=\omega_{*}=0$ also for
$|y|=1$, hence assertion~(a).

To prove assertion~(b), let us take $y\in B(0,1)$. Then, there
exists $\e\in (0,1/4)$ such that $y\in B(0,1-4\e)$. Since
$1-3\e>1-4\e$, there is $r_2(y)>0$ such that $B(y,r_2(y))\subset
B(0,1-3\e)$ and we deduce from \eqref{ld24} that there exists
$s_2(\e):=\log{(\tau_5(\e)+1)}>0$ such that $\omega(s,z)\ge\mu_\e$
for any $s\geq s_2(\e)$ and $z\in B(y,r_2(y))$. We now pick
sequences $(\sigma_n)_{n\geq1}$, $(\lambda_n)_{n\geq1}$ and
$(z_n)_{n\geq1}$ such that $\sigma_n\to0$, $\lambda_n\to\infty$,
$z_n\to y$, and $\omega(\sigma_n+\lambda_n,z_n)\to \omega_{*}(y)$.
Then there exists again $n_2(y)>0$ such that
$\sigma_n+\lambda_n>s_2(y)$ and $z_n\in B(y,r_2(y))$ for any $n\geq
n_2(y)$. Consequently $\omega(\sigma_n + \lambda_n,z_n) \geq \mu_\e$
for any $n\geq n_2(y)$. This readily implies that $\omega^{*}(y)\geq
\omega_{*}(y)\geq\mu_\e>0$, hence (b) is proved.

We follow the lines of \cite{LV07} and introduce
\begin{equation}
W_{*}(y)=\frac{p-1}{p-2} \omega_{*}(y)^{(p-2)/(p-1)}, \quad W^{*}(y)=\frac{p-1}{p-2} \omega^{*}(y)^{(p-2)/(p-1)},
\end{equation}
for any $y\in B(0,1)$. From Proposition~\ref{eik}, it follows that
$W_{*}$ and $W^{*}$ are respectively viscosity supersolution and
subsolution of the eikonal equation
$$
|\nabla \zeta|=1 \quad \mbox{in} \ B(0,1),
$$
with boundary conditions $W^{*}(y)=W_{*}(y)=0$ for $|y|=1$ and are both positive in $B(0,1)$. Using
the comparison principle of Ishii \cite{I87}, we find that
$W^{*}(y)\leq W_{*}(y)$, hence they should be equal by \eqref{c6}.
It follows that $\omega_{*}=\omega^{*}=W$ in $B(0,1)$, where $W$ is
the viscosity solution to \eqref{eqlim}
$$
|\nabla W|^{p-1} - W = 0 \quad\mbox{ in }\quad B(0,1)\,, \qquad W =0 \quad\mbox{ on }\quad \partial B(0,1)\,,
$$
which is actually explicit and given by
$$
W(x):=\left(\frac{p-2}{p-1}\ (1-|x|)_{+}\right)^{(p-1)/(p-2)},
$$
as stated in Theorem~\ref{asympt.main}. In addition, the equality
$\omega_*=\omega^{*}$ and \eqref{cvld21} entail the convergence of
$\omega(s)$ as $s\to \infty$ towards $W$ in $L^\infty(\real)$ by
Lemma~4.1 in \cite{Bl94} or Lemma~V.1.9 in \cite{BdCD97}. We end the
proof by rephrasing the two scaling steps and arriving in this way
to \eqref{main.asympt}. \qed

\subsection{Proof of Theorem~\ref{asympt.main}: $N\ge 2$}\label{subsec4.3}

We now prove Theorem~\ref{asympt.main} for radially non-increasing
initial data to the problem posed in dimension $N\ge 2$. We follow
the same steps as in dimension $N=1$, and we only indicate below the
main differences that appear. These differences are mainly given by
the appearance of the new term
\begin{equation}\label{newterm}
\frac{N-1}{r}|\partial_r w|^{p-2}\partial_r w \,, \quad r=|y|\,,
\end{equation}
in the radial form of the $p$-Laplacian term. As we shall see,
performing carefully the same steps as for dimension $N=1$, we find
that this term does not change anything in an essential way. We
follow the same division into steps as the case $N=1$.

\medskip

\noindent \textbf{Step 1.}  Thanks to the construction performed in Section~\ref{TWN}, this step is the same as in dimension $N=1$.

\medskip

\noindent \textbf{Step 2.} Due to the appearance of the extra term
\eqref{newterm} in the radial form of the equation \eqref{c1}, we
check by direct calculation that, in dimension $N\ge 2$, the
function $F_R$ given by formula \eqref{FR1} is now a strict
supersolution to \eqref{c1} in $Q$. Indeed, for $y\ne 0$,
$$
\partial_\tau F_R -\frac{1}{1+\tau}\left(\Delta_{p}F_R+y\cdot\nabla F_R-\frac{p-1}{p-2}F_R\right)+|\nabla F_R|^{p-1}-F_R =\frac{N-1}{(1+\tau) |y|} F_R\,.
$$
Moreover, its singularity at $y=0$ is now stronger. This seems to
introduce a new difficulty, but we will see that it can be handled
by the same perturbation techniques. Let us notice at this moment
that $F_R$ can be used for upper bounds in the same way as in the
case $N=1$, and that $F_R$ still solves the limit Hamilton-Jacobi
equation \eqref{limit}.

\medskip

\noindent \textbf{Step 3.} In order to construct subsolutions
starting from the family of functions $F_R$, we follow again the
ideas of the case $N=1$. The calculations will be different in some
points. We consider again the damped family $F_{R,\vartheta,\beta}$
defined in \eqref{marsupilami} for $R\in (0,1)$, $\vartheta\in
(0,1)$, and $\beta\in (1/2,1]$. For $y\neq0$ we have
\begin{equation*}
\begin{split}
Y&:=\partial_\tau F_{R,\vartheta,\beta}-\frac{1}{1+\tau}\left(\Delta_{p}F_{R,\vartheta,\beta}+y\cdot\nabla F_{R,\vartheta,\beta}-\frac{p-1}{p-2}F_{R,\vartheta,\beta}\right)+|\nabla F_{R,\vartheta,\beta}|^{p-1}-F_{R,\vartheta,\beta}\\
&= \vartheta F_{R,1,\beta}^{1/(p-1)} \left[\frac{\beta-\vartheta^{p-2}}{1+\tau} + \frac{(N-1)\vartheta^{p-2}}{(1+\tau) |y|} F_{R,1,\beta}^{(p-2)/(p-1)} - (1-\vartheta^{p-2}) F_{R,1,\beta}^{(p-2)/(p-1)} \right].
\end{split}
\end{equation*}
At this point, we further assume that $|y|>r_*/2$, the radius $r_*$ being defined in Corollary~\ref{cor:bfbsb}, and that
\begin{equation}\label{shortofidea}
\vartheta^{p-2} \le \frac{(1-\beta) r_*}{2(N-1)}\,.
\end{equation}
Since $F_{R,1,\beta}\le 1$, we obtain
\begin{equation*}
\begin{split}
Y & \le \vartheta F_{R,1,\beta}^{1/(p-1)} \left[\frac{\beta-\vartheta^{p-2}}{1+\tau} + \frac{2(N-1)\vartheta^{p-2}}{(1+\tau) r_*} - (1-\vartheta^{p-2}) F_{R,1,\beta}^{(p-2)/(p-1)} \right]\\
& \le \vartheta(1-\vartheta^{p-2}) F_{R,1,\beta}^{1/(p-1)} \left[\frac{1}{1+\tau}- \frac{p-2}{p-1}\left(\frac{\beta(\tau+R)}{\tau+1}-|y|\right)\right]\,,
\end{split}
\end{equation*}
from which we conclude that
\begin{equation}\label{stillshortofidea}
\begin{minipage}{10cm}
$F_{R,\vartheta,\beta}$ is a classical subsolution to \eqref{c1} in \\ $\{ (\tau,y)\in Q \ : \ \tau\ge\tau_2(R,\beta)\,, \ (r_*/2)<|y|\le K_{R,\beta}(\tau) \}$\,,
\end{minipage}
\end{equation}
where $\tau_2(R,\beta)$ and $K_{R,\beta}(\tau)$ are still given by
\eqref{pam}. We now proceed as in the one dimensional case to
establish \eqref{ld23} for all $\vartheta\in (0,1)$ satisfying
\eqref{ld22} along with
$$
\vartheta^{p-2} \le \frac{\e r_*}{N-1}\,,
$$
for \eqref{shortofidea} to be satisfied.

\medskip

\noindent \textbf{Steps 4 \& 5.} The final steps of the proof are similar to the one dimensional case. \qed


\section{Arbitrary initial data}\label{sec:aid}

So far, we have proved Theorems~\ref{asympt.main}
and~\ref{main.posit} for radially non-increasing initial data
satisfying \eqref{a3}. We now extend these two results to general
initial data satisfying \eqref{a3}.

\medskip

\noindent\textbf{Proof of Theorems~\ref{asympt.main}
and~\ref{main.posit}.} Since $u_0\not\equiv 0$, there are
$x_0\in\real^N$, $r_0>0$, and $\eta_0>0$ such that $u_0(x)\ge
2\eta_0$ for $x\in B(x_0,r_0)$. Then, there exists a radially
non-increasing initial condition $\tilde{u}_0$ satisfying \eqref{a3}
but with support in $B(0,r_0)$ and such that
$\|\tilde{u}_0\|_\infty\le\eta_0$ and $\tilde{u}_0(x) \le
u_0(x-x_0)$ for $x\in\real^N$. Similarly, there is a radially
non-increasing initial condition $\tilde{U}_0$ satisfying \eqref{a3}
but with support in $B(0,\tilde{R}_0)$ for some $\tilde{R}_0>R_0$
and such that $\tilde{U}_0(x)\ge \|u_0\|_\infty$ for $x\in
B(0,R_0)$. Denoting the solutions to \eqref{a1} by $\tilde{u}$ and
$\tilde{U}$ with initial conditions $\tilde{u}_0$ and $\tilde{U}_0$,
respectively, the comparison principle and the translational
invariance of \eqref{a1} ensure that
\begin{equation}\label{prunelle}
\tilde{u}(t,x+x_0) \le u(t,x) \le \tilde{U}(t,x)\,, \quad (t,x)\in Q\,.
\end{equation}
Moreover, since
$$
\left| \left( 1 - \frac{(p-2)|x+x_0|}{\log{t}}\right)_+^{(p-1)/(p-2)} - \left( 1 - \frac{(p-2)|x|}{\log{t}} \right)_+^{(p-1)/(p-2)} \right| \le \frac{(p-1)|x_0|}{\log{t}}\,,
$$
and Theorems~\ref{asympt.main} and~\ref{main.posit} apply to both
$\tilde{u}$ and $\tilde{U}$, the expected results follow from
\eqref{prunelle}. \qed

\section*{Appendix. Some results about viscosity solutions}

\setcounter{section}{7}

We state, for the sake of completeness, some standard results in the
theory of viscosity solutions, that we use in the proof of Theorem
\ref{main.asympt}. The first one concerns the ``viscosity'' limit of
a family of small perturbations and can be found in
\cite[Theorem~4.1]{Bl94}.
\begin{proposition}\label{pert}
Let $u_{\e}$ be a viscosity subsolution (resp. a viscosity
supersolution) of the equation
$$
H_{\e}(x,u_{\e},\nabla u_{\e},D^{2}u_{\e})=0 \quad \mbox{in} \
\real^N,
$$
where $H_{\e}$ is uniformly bounded in all variables and degenerate
elliptic. Suppose that $\{u_{\e}\}$ is a uniformly bounded family of
functions. Then
\begin{equation}
u^{*}(x):=\limsup\limits_{(y,\e)\to(x,0)}u_{\e}(y)
\end{equation}
is a subsolution of the equation
\begin{equation}
H_{*}(x,u,\nabla u,D^2u)=0,
\end{equation}
In the same way,
$$
u_{*}(x):=\liminf\limits_{(y,\e)\to(x,0)}u_{\e}(y)
$$
is a supersolution of $ H^{*}(x,u,\nabla u,D^2u)=0$. Here, $H_{*}$
and $H^{*}$ are constructed in the same way as \ $u_{*}$ and
$u^{*}$.
\end{proposition}
In other words, this result can be applied to asymptotically small perturbations of a known equation, as we do in Section~\ref{Sect.scaling2}.

We also use the following result:
\begin{proposition}\label{eik}
Let $u\in C(\Omega)$ be a viscosity solution of
\begin{equation}
H(x,u,\nabla u)=0 \quad \mbox{in} \ \Omega,
\end{equation}
where $\Omega\subset\real^N$ and $H$ is a continuous function. If $\Phi\in C^{1}(\real)$ is an increasing function, then $v=\Phi(u)$ is a viscosity solution of
\begin{equation}
H\left(x,\Phi^{-1}(v(x)),(\Phi^{-1})'(v(x))\nabla v(x)\right)=0.
\end{equation}
\end{proposition}
The same result holds true for subsolutions and supersolutions and
can be found in \cite{Bl94}. In particular, we use this result in
order to pass from the Hamilton-Jacobi equation $|\nabla
u|^{p-1}-u=0$ to the standard eikonal equation $|\nabla v|=1$.
Finally, we also use the (now standard) comparison principle for
viscosity subsolutions and supersolutions of the eikonal equation,
that can be found in \cite{I87}.

\bigskip

\textsc{Acknowledgements}. The first and the third author are
supported by the Spanish Projects MTM2005-08760-C02-01 and
MTM2008-06326-C02-01. JLV was partially supported by the ESF
Programme "Global and geometric aspects of nonlinear partial
differential equations". This work was initiated while the second
author enjoyed the support and hospitality of the Departamento de
Matem\'aticas of the Universidad Aut\'onoma de Madrid.

\bibliographystyle{plain}

\end{document}